\documentclass[pamm,a4paper,fleqn]{w-art}
\usepackage{times,cite,w-thm}
\usepackage[T1]{fontenc}
\usepackage[utf8]{inputenc}
\usepackage{graphicx}
\usepackage{amsmath,amsfonts,mathtools}
\usepackage{dsfont,xfrac,sidecap,caption}

\usepackage[ruled, vlined]{algorithm2e}

\newcommand{\nt}{\vec{n}} 

\newcommand{\vt}{\mathbf{v}} 

\newcommand{\ut}{\mathbf{u}}

\newcommand{\sigmat}{\boldsymbol{\sigma}}
\newcommand{\Sigmat}{\boldsymbol{\Sigma}}

\newcommand{\Ft}{\mathbf{F}}

\newcommand{\bs}[1]{\boldsymbol{#1}}

\newcommand{\FL}{{\cal F}}
\newcommand{\SO}{{\cal S}}
\newcommand{\IN}{\Gamma}

\begin{document}

\TitleLanguage[EN]
\title[The short title]{On temporal homogenization in the numerical simulation of atherosclerotic plaque growth}

\author{\firstname{Stefan} \lastname{Frei}\inst{1}%
} 
\address[\inst{1}]{\CountryCode[DE]Department of Mathematics \& Statistics, University of Konstanz, Germany (stefan.frei@uni-konstanz.de)}
\author{\firstname{Alexander} \lastname{Heinlein}\inst{2}%
}
\address[\inst{2}]{\CountryCode[NL]Delft Institute of Applied Mathematics, Delft University of Technology, The Netherlands (a.heinlein@tudelft.nl)}
\author{\firstname{Thomas} \lastname{Richter}\inst{3}%
}
\address[\inst{3}]{\CountryCode[DE]Institute for Analysis and Numerics, Otto-von-Guericke-University, Germany (thomas.richter@ovgu.de)}
\AbstractLanguage[EN]
\begin{abstract}
A temporal homogenization approach 
for the numerical simulation of atherosclerotic plaque growth is extended to fully coupled fluid-structure interaction (FSI) simulations. The numerical results indicate that the two-scale approach yields significantly different results compared to a simple heuristic averaging, where only stationary long-scale FSI problems are solved, confirming the importance of incorporating stress variations on small time-scales. In the homogenization approach, a periodic fine-scale problem, which is periodic with respect to the heart beat, has to be solved for each long-scale time step. Even if no exact initial conditions are available, periodicity can be achieved within only 2--3 heart beats by simple time-stepping.

\end{abstract}
\maketitle                   

\section{Introduction}

The growth of atherosclerotic plaque, leading to various cardiovascular diseases, happens on a long time scale of several months up to years. On the other hand, the growth of plaque is driven by the wall shear stress distribution, which varies significantly within each heartbeat, that is, over each second. It is unfeasible to perform numerical simulations ranging over multiple months with time step sizes in the order of magnitude of milliseconds, which makes temporal homogenization techniques indispensable. Here, the temporal homogenization approach from~\cite{FreiRichter2020} is extended to fully coupled fluid-structure interaction (FSI) simulations. In this approach, for each long-scale time step, a fine-scale FSI solution, which is periodic with respect to the heart beat, has to be computed. In order to justify this approach and confirm its efficiency, we discuss two important questions. On the one hand, we investigate whether long-scale simulations benefit from including averaged fine-scale information. On the other hand, we investigate the periodicity of the fine-scale problem. In particular, the initial data of a periodic solution on the short-scale is unknown a priori. Due to the diffusivity of the flow problem a remedy is a simple time-stepping approach from an approximate initial state. We investigate how many heart beats are necessary to obtain an almost-periodic state which can be used as initial data.

\section{Fluid-structure interaction}

We consider a partition of an overall domain $\Omega(t) = \FL(t) \cup \IN(t)\cup \SO(t)$ into a fluid part $\FL(t)$, an interface $\IN(t)$ and a solid part $\SO(t)$.
Blood flow through a vessel wall is described by the following FSI system:
\begin{equation} \label{fullplaquemodel}
  \begin{aligned}
      \rho_f(\partial_t \vt_f+\vt_f\cdot\nabla\vt_f) -
      \operatorname{div}\,\sigmat_f&= 0, &\qquad
      \operatorname{div}\,\vt_f&=0 \qquad\text{ in }\FL(t), \\
      \rho_s (\partial_t\hat{\vt}_s) -
      \operatorname{div}\,(\hat{\Ft}_e\hat{\Sigmat}_e) &=0, &\qquad
      \partial_t \hat\ut_s - \hat\vt_s &=0 \qquad \text{ in }\hat{S}, \\
      \sigmat_f\nt_f+\sigmat_s\nt_s &=0, &\qquad
      \vt_f&=\vt_s \quad\;\; \text{ on }\IN(t).
  \end{aligned}
\end{equation}
Here, $\vt_f$ and $\hat{\vt}_s$ stand for the fluid and solid velocity and $\hat{\ut}_s$ for the solid displacement. Quantities with a ``hat'' are defined in Lagrangian coordinates, while quantities without a "hat" are given in the current Eulerian coordinate framework. They correspond to each other by a $C^{1,1}$-diffeomorphism $\xi: \hat{\Omega}\to \Omega(t)$ and the relation $\hat{f} = f\circ \xi$, which also defines the solid deformation gradient $\hat{F}_s = \nabla \xi|_{\hat{S}} = I + \nabla \hat{\ut}_s$. By $\rho_f$ and $\rho_s$ we denote the densities of blood and vessel wall and by $\nt_f$ and $\nt_s$ the outer normals of the fluid and solid domain, respectively.

\begin{figure}
\centering
\begin{picture}(0,0)%
\includegraphics{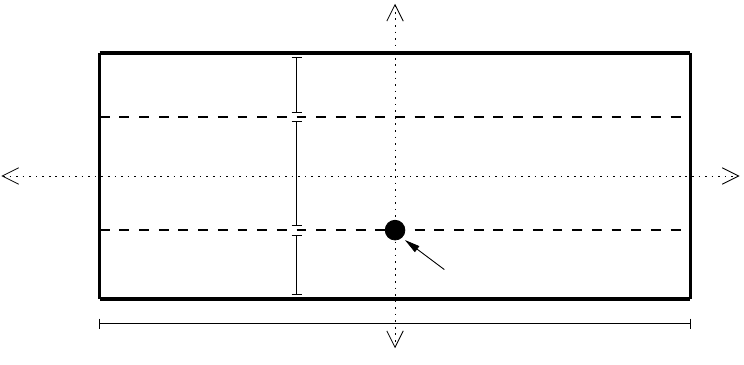}%
\end{picture}%
\setlength{\unitlength}{2072sp}%
\begingroup\makeatletter\ifx\SetFigFont\undefined%
\gdef\SetFigFont#1#2{%
  \fontsize{#1}{#2pt}%
  \selectfont}%
\fi\endgroup%
\begin{picture}(6774,3312)(439,-2911)
\put(991,-421){\makebox(0,0)[lb]{\smash{{\SetFigFont{9}{10.8}{\color[rgb]{0,0,0}$\Gamma_s$}%
}}}}
\put(991,-2221){\makebox(0,0)[lb]{\smash{{\SetFigFont{9}{10.8}{\color[rgb]{0,0,0}$\Gamma_s$}%
}}}}
\put(5311,-421){\makebox(0,0)[lb]{\smash{{\SetFigFont{9}{10.8}{\color[rgb]{0,0,0}$\SO$}%
}}}}
\put(5311,-2221){\makebox(0,0)[lb]{\smash{{\SetFigFont{9}{10.8}{\color[rgb]{0,0,0}$\SO$}%
}}}}
\put(6841,-421){\makebox(0,0)[lb]{\smash{{\SetFigFont{9}{10.8}{\color[rgb]{0,0,0}$\Gamma_s$}%
}}}}
\put(6841,-2221){\makebox(0,0)[lb]{\smash{{\SetFigFont{9}{10.8}{\color[rgb]{0,0,0}$\Gamma_s$}%
}}}}
\put(1396,-2896){\makebox(0,0)[lb]{\smash{{\SetFigFont{9}{10.8}{\color[rgb]{0,0,0}10}%
}}}}
\put(3241,-421){\makebox(0,0)[lb]{\smash{{\SetFigFont{9}{10.8}{\color[rgb]{0,0,0}1}%
}}}}
\put(3241,-2221){\makebox(0,0)[lb]{\smash{{\SetFigFont{9}{10.8}{\color[rgb]{0,0,0}1}%
}}}}
\put(901,-1051){\makebox(0,0)[lb]{\smash{{\SetFigFont{9}{10.8}{\color[rgb]{0,0,0}$\Gamma_\text{in}$}%
}}}}
\put(6841,-916){\makebox(0,0)[lb]{\smash{{\SetFigFont{9}{10.8}{\color[rgb]{0,0,0}$\Gamma_\text{out}$}%
}}}}
\put(3196,-1456){\makebox(0,0)[lb]{\smash{{\SetFigFont{9}{10.8}{\color[rgb]{0,0,0}2}%
}}}}
\put(5491,-1051){\makebox(0,0)[lb]{\smash{{\SetFigFont{9}{10.8}{\color[rgb]{0,0,0}$\FL$}%
}}}}
\put(4186,-601){\makebox(0,0)[lb]{\smash{{\SetFigFont{9}{10.8}{\color[rgb]{0,0,0}$\IN$}%
}}}}
\put(4186,-1636){\makebox(0,0)[lb]{\smash{{\SetFigFont{9}{10.8}{\color[rgb]{0,0,0}$\IN$}%
}}}}
\put(4501,-2131){\makebox(0,0)[lb]{\smash{{\SetFigFont{6}{7.2}{\color[rgb]{0,0,0}$(0,-1)$}%
}}}}
\put(4096,-1096){\makebox(0,0)[lb]{\smash{{\SetFigFont{6}{7.2}{\color[rgb]{0,0,0}$(0,0)$}%
}}}}
\end{picture}%
\caption{\label{fig:conf} Sketch of the simulation domain: a two-dimensional channel where the plaque growth is centered at the black dot at $(0,-1)$.}
\end{figure}
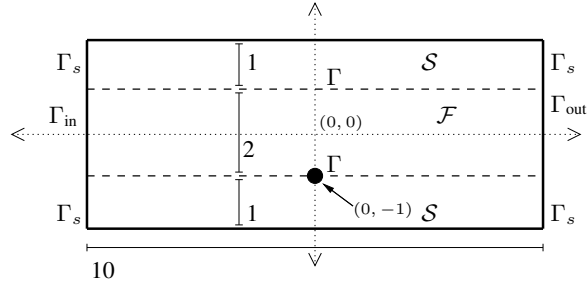

The Cauchy stress tensors $\sigmat_f$ and $\sigmat_s$ for the fluid and solid, respectively, are given by the respective constitutive models. Here, we employ the Saint Venant--Kirchhoff model with Lam\'e parameters $\mu_s$ and $\lambda_s$ combined with a growth model to describe the solid; see Section~\ref{sec:solid}. The blood is modeled as an incompressible Newtonian fluid with Cauchy stresses
\begin{equation}\label{fluidstresses-euler}
  \sigmat_f=\rho_f\nu_f(\nabla\vt_f+\nabla\vt_f^T)-p_f I,
\end{equation}
with the kinematic viscosity $\nu_f$, which results in the Navier--Stokes equations.
A sketch of the domain of interest is given in Figure~\ref{fig:conf}. The boundary data of our FSI problem is given by 
\begin{align}\label{problem:long-scale-boundary}
	\vt_f = \vt^\text{in} \text{ on }\Gamma^\text{in}_f,
	\quad
	\rho_f\nu_f\nt\cdot\nabla\vt_f-p_f\nt =0
	& \text{ on }\Gamma^\text{out}_f,\quad
	\hat{\ut}_s=0\text{ on }\Gamma_s,
\end{align}
where $\nt$ is the outward facing normal vector and the inflow $\vt^\text{in}$ is specified below.
  
\section{Modeling of solid growth} \label{sec:solid}

The modeling of plaque growth at the vessel wall is a complex task that involves the interaction of many different molecules and species; see for example~\cite{SilvaJaegerNeussRaduSequeira2020, SonnerPhD}. In this contribution, our focus does not lie on realistic modeling of plaque growth but on methodological questions of the temporal homogenization approach. Hence, we employ a greatly simplified growth model. Specifically, we concentrate on the influence of the concentration of foam cells $c_s$ on the growth. The rate of formation of these cells depends on the wall shear stress $\sigmat_f^{WS}$ at the vessel wall. In our simplified model, this is included by means of a simple ordinary differential equation (ODE) for $c_s$,  %
\begin{equation}\label{reaction}
	\partial_t c_s = \gamma(\sigmat_f^{WS},c_s)\coloneqq {\alpha}
	\big(1+c_s\big)^{-1}\big(1+|\sigmat_f^{WS}|^2\big)^{-1},\quad
	\sigmat_f^{WS}\coloneqq \sigma_0^{-1}\int_\Gamma\rho\nu
	\big(I_d-\nt_f\nt_f^T\big)(\nabla\vt+\nabla\vt^T)\nt_f\,\text{d}o,
\end{equation}
where $\sigma_0$ and {periodic $\alpha$} are parameters of the growth model; see also~\cite{FreiRichter2020}. 

We model solid growth by a multiplicative splitting of the deformation gradient $\hat \Ft_s$ into an elastic part $\hat\Ft_e$ and a growth function $\hat\Ft_g$~\cite{RodriguezHogerMcCulloch1994,YangJaegerNeussRaduRichter2015,FreiRichterWick2016}
 \begin{equation}\label{elasticgrowth-reference}
  \hat \Ft_s=\hat\Ft_e\hat\Ft_g\quad\Leftrightarrow\quad
  \hat\Ft_e = \hat\Ft_s\hat\Ft_g^{-1} =
  [I+\hat\nabla\hat\ut_s]\hat\Ft_g^{-1}.
\end{equation}
In our simplified model configuration, we consider the following growth function depending on $c_s$
\begin{equation}\label{num:1:growth}
  \hat g(\hat x,\hat y,t) =  1+ c_s(t) \exp\left(-\hat
  x^2\right)(2-|\hat y|),\quad 
  \hat\Ft_g(\hat x,\hat y,t):=\hat g(\hat x,\hat y,t)\,I.
\end{equation}
This means that the shape and position of the plaque growth is prescribed, but the growth rate depends on the variable $c_s$. As the simulation domain is centered around the origin, growth takes place in the middle of the domain; see Figure~\ref{fig:conf}.

From~\eqref{elasticgrowth-reference} and~\eqref{num:1:growth}, it follows that $\hat\Ft_e:=\hat g^{-1} \hat\Ft_s$. The Piola--Kirchhoff stresses of the St.\,Venant-Kirchhoff model take the form \\[-0.5cm]
\begin{align*}\label{Piola}
	\hat\Ft_e\hat\Sigmat_e = 2\mu_s \hat g^{-1}\hat{\bs{F}}_s\hat{\bs{E}}_e + \lambda_s  \hat g^{-1}\operatorname{tr}(\hat{\bs{E}}_e)\hat{\bs{F}}_s,
	\qquad
	\hat{\bs{E}}_e = \frac{1}{2}(\hat g^{-2}\hat
      {\bs{F}}_s^T\hat{\bs{F}}_s-I) 
\end{align*}
which correspond to the Cauchy stresses by
$ 
 \sigmat_s(x) = \hat{\sigmat}_s(\hat{x}) = \hat{J}_e^{-1} \hat{\Ft}_e\hat{\Sigmat}_e(\hat{x}) \hat{\Ft}_e^T
$.

\section{Numerical framework}

We use an arbitrary Lagrangian--Eulerian (ALE) approach to solve the FSI problem~\eqref{fullplaquemodel}, and make use of the symmetry of the configuration to restrict the simulation to the lower half of the computational domain. For the details, we refer to~\cite{FreiRichterWick2016}.

A resolution of the short-scale dynamics with a scale of milliseconds to seconds over the complete time interval of interest $[0,T]$ with $T$ being several months up to a year is unfeasible, even for the simplified two-dimensional configuration considered here. For $T=200$ days and a relatively coarse short-scale time step of $\delta \tau = 0.02\,$s, a total of $8.64\cdot 10^8$ time steps would be required, where each step corresponds to the solution of a mechano-chemical FSI problem.

This dilemma is frequently solved by considering a heuristic averaging approach: Assuming that a stationary limit of the FSI system on the short scale is reached, a stationary FSI problem can be solved on a longer scale (e.g., $\delta t\approx 1$ day). The wall-shear stress $\overline{\sigmat}_f^{WS}$ of the stationary FSI problem is then used to advance the foam cell concentration based on~\eqref{reaction}.
  
Note that $\overline{\sigmat}_f^{WS}$ is not necessarily a good approximation of $\sigmat_f^{WS}$, which depends on the pulsating blood flow, as shown in \cite{FreiRichterWick2016} and analyzed in~\cite{FreiRichter2020, SonnerPhD}. The stationary {FSI} model reads  
\begin{equation}\label{statmodel}
\begin{aligned}
	\rho_f\bar{\vt}_f\cdot\nabla\bar{\vt}_f -
    \operatorname{div}\,\bar{\sigmat}_f&= 0, \quad
    &\operatorname{div}\,\bar{\vt}_f&=0 \quad\text{ in }\FL(t),&\quad
    \operatorname{div}\,(\bar{\Ft}_e\bar{\Sigmat}_e) &=0 \quad\text{ in }\hat{S}\\
    \bar{\sigmat}_f\nt_f+\bar{\sigmat}_s\nt_s &=0, \qquad &\bar{\vt}_f&=\bar{\vt}_s \;\; \text{ on }\IN(t), &\quad
    \partial_t \bar{c}_s &= \gamma(\overline{\sigmat}_f^{WS}, \bar{c}_s).
\end{aligned}  
\end{equation}
The corresponding algorithm is given as Algorithm~\ref{averaging}.

\begin{algorithm}[t]
\caption{Simple Heuristic Averaging\label{averaging}}
  	Initialize $c_s^0=0$. Set time step size $\delta t = 1\,\text{day}=86\,400\,s$. \\
  	\For{$n=1,2,\dots$}{ 
  	\begin{enumerate} \setlength{\itemsep}{-1.5ex}
  	\item Solve quasi-stationary Long-scale problem~\eqref{statmodel}:
    	\vspace{-1ex}\[
    	\{\bar{c}_s^{n-1}\} \mapsto 
    	\{\bar{\vt}^{n},\bar{\ut}^{n},\bar{p}^{n}\}
    	\]
  	\item Compute the wall stress in main stream direction
    	\vspace{-1ex}\begin{equation*}
    	\overline{\sigmat}_f^{WS,n} = \int_{\Gamma_i}
    	\sigmat_f(\bar{\vt}^n,\bar{p}^n)\nt_f\cdot\vec{e}_1\,\text{d}o\label{WSS}
    	\end{equation*}
  	\item Update the foam cell concentration
    	\vspace{-1ex}\[
    	\bar{c}_s^{n} = \bar{c}_s^{n-1} + \delta t \gamma(\overline{\sigmat}_f^{WS,n}, \bar{c}_s^{n-1})
    	\]\vspace{-1cm}
    \end{enumerate}}
\end{algorithm}

A more accurate two-scale approach has been presented by Frei and Richter in~\cite{FreiRichter2020}, where a periodic-in-time short-scale problem needs to be solved in each time step of the long (\textit{macro}) scale, e.g., each day. The growth function $\gamma(\sigmat_f^{WS})$ is then averaged by integrating over one period of the heart beat, and the average is denoted by $\overline{\gamma}(\sigmat_f^{WS})$.

A difficulty lies in the solution of the periodic short-scale problem. In general, the exact initial data of the periodic solution is not know a priori. Instead, if a reasonable guess for starting values $\vt^0$ and $\ut_s^0$ is available on the short-scale, the FSI system may convergence to a periodic state. After each cycle, we can check if the solution is sufficiently close to a periodic state. In this work, we apply a stopping criterion based on the computed averaged growth value (where $k$ denotes the iteration index)
\begin{equation}
	|\overline{\gamma}(\sigmat_f^{WS,k}) - \overline{\gamma}(\sigmat_f^{WS,k-1})| < \epsilon_p.
\end{equation}
  
\begin{algorithm}
\caption{Two-Scale Algorithm\label{twoscale}}
  	Initialize $c_s^0=0$.  Set time step size $\delta t = 1\,\text{day}=86\,400\,\text{s}$. \\
  	\For{$n=1,2,\dots$}{ 
  	\begin{enumerate} \setlength{\itemsep}{-1.5ex}
  	\item[1.)] Solve quasi-stationary Long-Scale Problem~\eqref{statmodel}:
    	\vspace{-1ex}\[
    	\{c_s^{n-1}\} \mapsto 
	    \{\vt^{n},\ut^{n},p^{n}\}
    	\]
    \item[2.)] Set suitable starting values $\vt^{0,0}$ and $\ut^{0,0}$.\\
    	\While{$|\overline{\gamma}(\sigmat_f^{WS,k,m}) - \overline{\gamma}(\sigmat_f^{WS,k-1,m})| > \epsilon_p$}{ 
	    \begin{itemize} \setlength{\itemsep}{-1.5ex}
  	\item[2.a)]
  		Solve Short-Scale Problem~\eqref{fullplaquemodel} 
    	in $I_n =(\tau_{n-1},\tau_{n-1}+1\,\text{s})$
    	\vspace{-1ex}\[    
	    \{\vt^{k,0},\ut^{k,0},c_s^{n-1}\}\mapsto 
    	\{\vt^{k,m},\ut^{k,m},p^{k,m}\},\; m=1,\dots,N_s
	    \]
  	\item[2.b)] Compute average growth function
    	\vspace{-1ex}\[
    	\overline{\gamma}(\sigmat_f^{WS}, c_s^{n-1}) = \frac{1}{N_s}\sum\nolimits_{m=1}^{N_s}
	    \gamma(\sigmat_f^{WS,m}, c_s^{n-1}), \quad k \gets k+1
    	\]\vspace{-0.5cm}
		\end{itemize}}
  	\item[3.)] Update the foam cell concentration
    	\vspace{-1ex}\[
	    c_s^{n} = c_s^{n-1} + \delta t\,\overline{\gamma}(\sigmat_f^{WS,n}, c_s^{n-1})
    	\]\vspace{-1cm}
    \end{enumerate}
}
\end{algorithm}

The algorithm is summarized in Algorithm~\ref{averaging}. We consider different possibilities to set the initial values $\vt^{0,0}$ and $\ut^{0,0}$ in step 2 of Algorithm 2. Either the variables $\vt^{k,m}$ and $\ut^{k,m}$ from the previous day can be used (``micro'' strategy), or the values $\vt^{n}$, $\ut^{n}$ from the previous long-scale step (``macro'' strategy). In the latter case, we set $\vt_f^{n-1}$ to zero in order to match the inflow boundary conditions at the begin of the period. If the ``micro'' strategy is used, step 1 of Algorithm~\ref{twoscale} can be skipped.

\section{Numerical example}

We consider a channel with length $10\,$cm and an initial width $\omega(0)$ of $2\,$cm as illustrated in Figure~\ref{fig:conf}. The solid walls on top and bottom have an initial thickness of 1 cm each. Fluid density and viscosity are given by $\rho_f=1\,\text{g}/\text{cm}^3$ and $\nu_f=0.04\,\text{cm}^2/\text{s}$, respectively, the solid density is $\rho_s=1\,\text{g}/\text{cm}^3$, and the Lam\'e parameters are $\mu_s=10^4$ and $\lambda_s=4\cdot 10^4\,\text{dyn}/\text{cm}^2$. The growth parameters are $\sigma_0=30 \frac{\text{g\,cm}}{\text{s}^2}$ and {periodic $\alpha = 0.0432 \frac{1}{\rm day} = 5 \cdot 10^{-7} \frac{1}{\rm s}$}. We prescribe a pulsating velocity inflow profile on $\Gamma_f^\text{in}$ by
\begin{equation}\label{num:1:inflow}
	\vt_f^\text{in}(t,x,y) = 30\begin{pmatrix}
   	\sin(\pi t)^2 (1-y^2)\\ 0
  	\end{pmatrix}  \text{cm}/\text{s}.
\end{equation}

We use $Q_2$ equal-order finite elements with local projection stabilization (LPS) stabilization~\cite{BeckerBraack2001} for all variables on a mesh consisting of 160 rectangular grid cells, which corresponds to a total of $3\,157$ degrees of freedom. The time step sizes are chosen as $\delta \tau=0.02\,$s and $\delta t=1\,$day, the tolerance for periodicity $\epsilon_p=10^{-3}$. All the computational results have been obtained with the finite element library Gascoigne3d~\cite{Gascoigne}.

\begin{figure}
\includegraphics[width=0.45\textwidth]{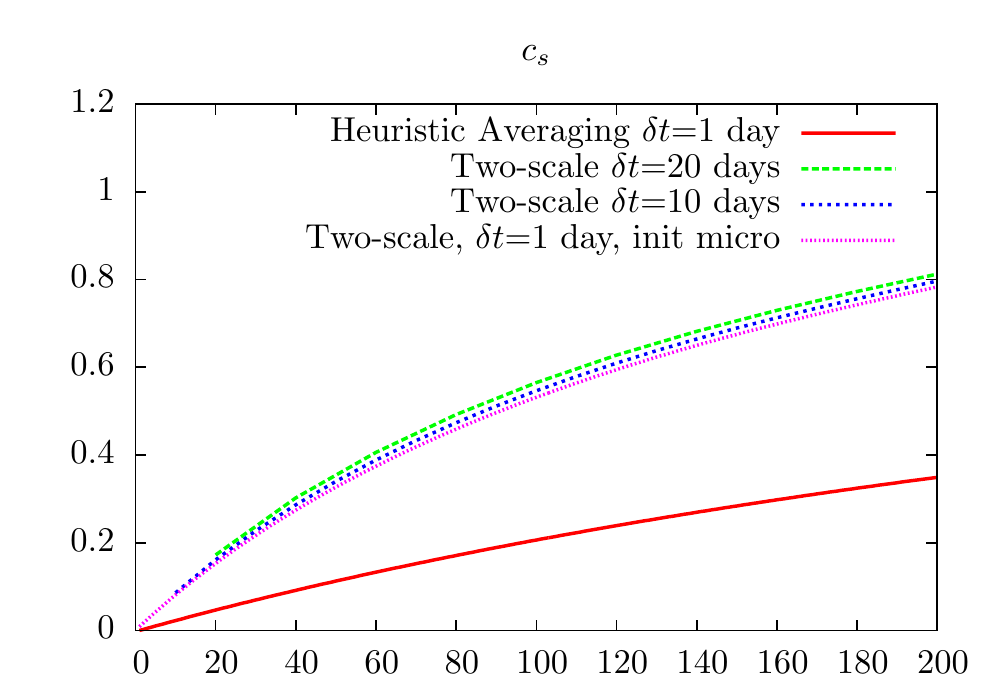}
\includegraphics[width=0.45\textwidth]{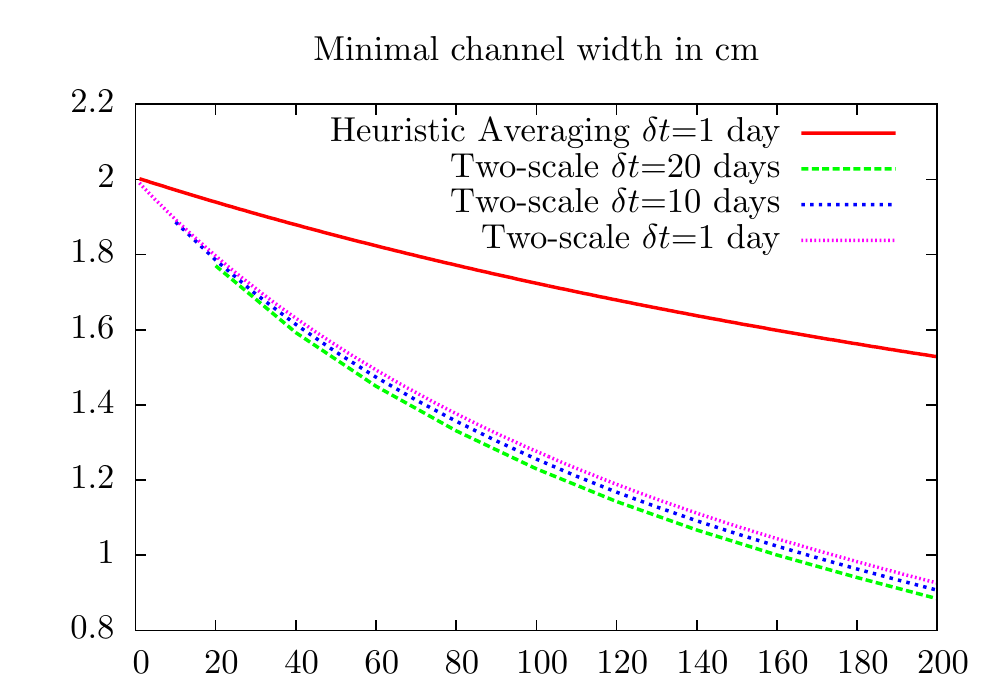}\\
\includegraphics[width=0.45\textwidth]{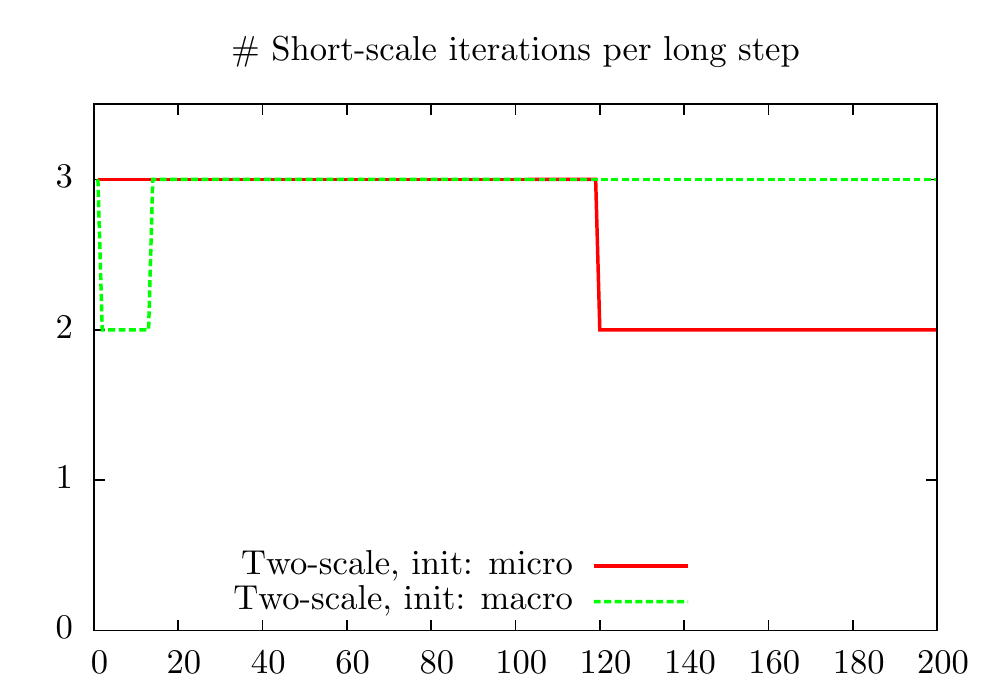}
\includegraphics[width=0.45\textwidth]{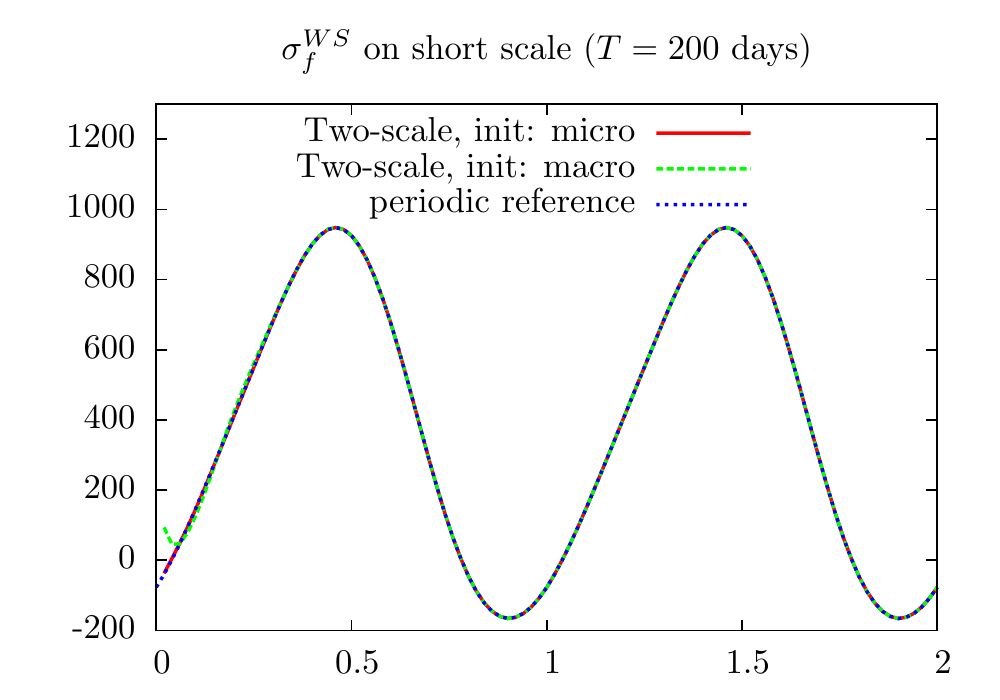}
\caption{\label{fig:NumEx} Comparison of the simple heuristic averaging approach (Algorithm 1) and the two-scale approach with different initialization strategies and long-scale time step sizes $\delta t$. \textbf{Top left}: growth function over time. \textbf{Top right}: channel width over time. \textbf{Bottom left}: numbers of short-scale iterations per time step to reach quasi-periodicity. \textbf{Bottom right}: wall shear stress on the short scale at $T=200$ days.}
\end{figure}

Our results are given in Figure~\ref{fig:NumEx}. We see that the pure averaging approach (Algorithm 1) underestimates the growth significantly. A very coarse discretization of the long time interval gives a much better approximation. Concerning the different initialization strategies, we observe small deviations in the wall shear stress in the very beginning, in particular for the ``macro'' (last long-step) strategy. Convergence to the periodic state, is however reached very quickly, see Figure~\ref{fig:NumEx}, bottom right. Only 2--3 cycles are required until the periodicity criterion is met. Figure~\ref{fig:NumEx}, bottom left, shows an advantage of the ``micro'' initialization strategy, where more frequently the stopping criterion is fulfilled after only two cycles.

\begin{acknowledgement}
\end{acknowledgement}

\bibliographystyle{plain}

\begin{thebibliography}{1}

\bibitem{BeckerBraack2001}
R~Becker and M~Braack.
\newblock A finite element pressure gradient stabilization for the {S}tokes
  equations based on local projections.
\newblock {\em Calcolo}, 38(4):173--199, 2001.

\bibitem{Gascoigne}
R.~Becker, M.~Braack, D.~Meidner, T.~Richter, and B.~Vexler.
\newblock The finite element toolkit {G}ascoigne3d.
\newblock http://www.gascoigne.de.

\bibitem{FreiRichter2020}
S.~Frei and T.~Richter.
\newblock Efficient approximation of flow problems with multiple scales in
  time.
\newblock {\em SIAM Multiscale Model Simul}, 18(2):942--969, 2020.

\bibitem{FreiRichterWick2016}
S.~Frei, T.~Richter, and T.~Wick.
\newblock Long-term simulation of large deformation, mechano-chemical
  fluid-structure interactions in {ALE} and fully {E}ulerian coordinates.
\newblock {\em J Comput Phys}, 321:874 -- 891, 2016.

\bibitem{RodriguezHogerMcCulloch1994}
E.K. Rodriguez, A.~Hoger, and A.D. McCulloch.
\newblock Stress-dependent finite growth in soft elastic tissues.
\newblock {\em J. Biomech.}, 4:455--467, 1994.

\bibitem{SilvaJaegerNeussRaduSequeira2020}
T.~Silva, W.~J{\"a}ger, M.~Neuss-Radu, and A.~Sequeira.
\newblock Modeling of the early stage of atherosclerosis with emphasis on the
  regulation of the endothelial permeability.
\newblock {\em J Theor Biol}, 496:110229, 2020.

\bibitem{SonnerPhD}
F.~Sonner.
\newblock {\em Temporal Multiscale Methods for a Model of Atherosclerosis}.
\newblock PhD thesis, FAU Erlangen-N{\"u}rnberg, 2021.

\bibitem{YangJaegerNeussRaduRichter2015}
Y.~Yang, W.~J{\"a}ger, M.~Neuss-Radu, and T.~Richter.
\newblock Mathematical modeling and simulation of the evolution of plaques in
  blood vessels.
\newblock {\em J Math Biol}, pages 1--24, 2014.

\end{thebibliography}

\end{document}